\documentclass[a4paper,11pt]{article}
\usepackage{graphicx}
\usepackage{amssymb}
\usepackage{amsthm}
\usepackage{subcaption}
\usepackage{amsfonts}
\usepackage{amsmath}
\usepackage{hyperref}
\usepackage{xcolor}
\usepackage{seqsplit}
\newcommand{\hexseq}[1]{{\footnotesize\ttfamily\seqsplit{#1}}}

\title{Searching for binary orientable sequences}
\author{Chris J. Mitchell\\Information Security Group, Royal Holloway, University of London\\
\url{www.chrismitchell.net}}
\date{24th October 2026}
\theoremstyle{plain}

\theoremstyle{definition}
\newtheorem{definition}{Definition}[section]

\newtheorem{example}{Example}[section]

\theoremstyle{remark}
\newtheorem{remark}{Remark}[section]

\begin{document}

\maketitle

\begin{abstract}
An aperiodic binary orientable sequence is a finite sequence of zeros and ones for which, for some
integer $n$, all the subsequences of $n$ consecutive symbols, read either forwards or backwards,
are distinct. Analogously, a periodic binary orientable sequence is an infinite periodic sequence
of zeros and ones for which, for some integer $n$, all the subsequences of $n$ consecutive symbols,
read either forwards or backwards, are distinct within a single period of the sequence. Such
sequences have potential application in position-location applications. For applications, the
longer the sequence (or the greater the period in the periodic case) the greater its utility, but
yet while constructions are known for such sequences, and upper bounds on their length (or period)
exist, the precise maximum length (or period) was previously unknown for anything other than small
values of $n$. In this note the results of LLM-conducted searches for the longest possible
aperiodic and periodic binary orientable sequences for $n\leq16$ are reported. In each case an
example of a sequence of the maximum length is also provided, and the maximum possible length is
now determined for every $n\leq16$ for both the periodic and aperiodic cases.
\end{abstract}

\section{Introduction} \label{section:introduction}

An aperiodic binary orientable sequence is a finite sequence of zeros and ones for which, for some
integer $n$ (the \emph{span}), all the subsequences of $n$ consecutive symbols, read in either
direction, are distinct. Analogously, a periodic binary orientable sequence is an infinite periodic
sequence of zeros and ones for which, for some integer $n$ (the \emph{span}), all the subsequences
of $n$ consecutive symbols, read either forwards or backwards, are distinct within a single period
of the sequence.

Such sequences were introduced in the early 1990s \cite{Burns92,Burns93,Dai93} in the context of
position-location applications.  Suppose such a sequence is encoded in some way in a surface, and
suppose that a device capable of reading $n$ consecutive symbols of the encoding wishes to learn
where it is and in which direction it is traversing the sequence. Then, because no $n$-tuple
repeats in either direction, reading $n$ consecutive symbols can be used to work out the location
in the sequence and the direction in which it is being traversed.

The searches described in this note were performed with the aid of the Claude
LLM\footnote{\url{https://claude.ai}}.

\section{Definitions}

We first give a formal definition of the sequences forming the focus of this note.  Note that
from here on we do not explicitly say that we are only dealing with the binary case.

\begin{definition}
Suppose $n\geq2$.  A finite sequence of length $\ell$, $(a_0,a_1,\dots,a_{\ell-1})$, where
$a_i\in\{0,1\}$ for every $i$, is said to be an \emph{aperiodic orientable sequence} of span $n$ if
all the $2(\ell-n+1)$ $n$-tuples $(a_i,a_{i+1},\dots,a_{i+n-1}$ and
$(a_{i+n-1},a_{i+n-2},\dots,a_i)$ ($0\leq i\leq \ell-n$) are distinct.
\end{definition}

\begin{example}
A trivial example for $n=2$ is $(01)$.  For $n=3$, a simple example is $(0011)$.
\end{example}

A periodic binary orientable sequence can be defined as follows.

\begin{definition}
Suppose $n\geq2$.  An infinite periodic sequence $(a_i)$ ($i\geq0$) of period $p$, where
$a_i\in\{0,1\}$ for every $i$, is said to be a \emph{periodic orientable sequence} of span $n$ if
all the $2p$ $n$-tuples $(a_i,a_{i+1},\dots,a_{i+n-1}$ and $(a_{i+n-1},a_{i+n-2},\dots,a_i)$
($0\leq i\leq p-1$) are distinct.
\end{definition}

\section{Previous results}

\subsection{The 1993 table}

The 1993 Burns and Mitchell paper \cite[Table 3]{Burns93} tabulates the results of searches for
aperiodic orientable sequences, and is reproduced below.  Table~\ref{table:burnstab} lists the
lengths of the longest known aperiodic orientable sequence for $4 \leq n \leq 16$.  For $n$ in
the range $4 \leq n \leq 7$ the length given is the length of the longest such sequence, since
for these values an exhaustive search was completed.

\begin{table}[htb]
\caption{Existence of aperiodic orientable sequences}  \label{table:burnstab}
\begin{center}
\begin{tabular}{cc} \hline
Span ($n$) & Sequence length ($\ell$) \\ \hline
                  4     &      8 \\
                  5     &      14 \\
                  6     &      26 \\
                  7     &      48 \\
                  8     &      108 \\
                  9     &      210 \\
                  10    &      440 \\
                  11    &      872 \\
                  12    &      1860 \\
                  13    &      3710 \\
                  14    &      7400 \\
                  15    &      15467 \\
                  16    &      31766 \\ \hline
\end{tabular}
\end{center}
\end{table}

\subsection{A recursive construction}

A 2022 paper \cite{Mitchell22} describes a simple recursive approach for generating aperiodic
orientable sequences.  This approach can be used to generate sequences with lengths given in
Table~\ref{table:AOS-family-table}.

\begin{table}[htb]
\caption{A family of aperiodic orientable sequences}  \label{table:AOS-family-table}
\begin{center}
\begin{tabular}{cc} \hline
Span ($n$) & Sequence length ($\ell$) \\ \hline
                  2     &      2 \\
                  3     &      4 \\
                  4     &      8 \\
                  5     &      14 \\
                  6     &      26 \\
                  7     &      48 \\
                  8     &      92 \\
                  9     &      178 \\
                  10    &      350 \\
                  $\cdots$ & $\cdots$ \\
                  $2r$ & $2^{2r-2}+(2^{2r-2}-1)/3+2r-1$ \\
                  $2r+1$ & $2^{2r-1}+(2^{2r-1}-2)/3+2r$ \\ \hline
\end{tabular}
\end{center}
\end{table}

It should be clear by comparing with Table~\ref{table:burnstab} that this method generates maximal
length sequences for $n\leq7$, but not for larger $n$.

\subsection{Constructions for periodic sequences}

Most of the prior work in this general area has focussed on the \emph{periodic} case, and the two
cases are closely related as the following remark shows.

\begin{remark}
Any $p+n-1$ consecutive symbols of a periodic orientable sequence form an aperiodic
orientable sequence.  That is, any method of construction for a periodic orientable sequence
automatically yields aperiodic sequences.
\end{remark}

A number of recent papers have considered construction methods, and period bounds, for periodic
sequences, including \cite{Gabric24b,Mitchell22}.  Table~\ref{table:Gabric} lists the periods of
the orientable sequences constructed by Gabri\'{c} and Sawada~\cite{Gabric24b}, together with the
length of the aperiodic orientable sequence which can be derived from it.

\begin{table}[htb]
\caption{The Gabri\'{c}-Sawada sequences}  \label{table:Gabric}
\begin{center}
\begin{tabular}{ccc} \hline
Span ($n$) & Period ($p$) & Aperiodic sequence length ($\ell=p+n-1$) \\ \hline
    5     &      6       &  10     \\
    6     &      16      &  21     \\
    7     &      36      &  42     \\
    8     &      92      &  99     \\
    9     &      174     &  182    \\
    10    &      416     &  425    \\
    11    &      844     &  854    \\
    12    &      1844    &  1855   \\
    13    &      3700    &  3712   \\
    14    &      7694    &  7707   \\
    15    &      15394   &  15408  \\
    16    &      31483   &  31498  \\
    17    &      63135   &  63151  \\
    18    &      128639  &  128656 \\
    19    &      257272  &  257290 \\
    20    &      519160  &  519179 \\  \hline
\end{tabular}
\end{center}
\end{table}

For further information on orientable sequences and their generation see Sawada's excellent web
page: \url{https://debruijnsequence.org/db/orientable}

\section{Searches}

\subsection{Aperiodic sequences}  \label{subsection:aperiodic_searches}

Claude used two methods of search: exhaustive search (for $n\leq7$) and integer programming ($8\leq
n\leq 16$), as follows.
\begin{itemize}
\item \emph{Exhaustive search}.  This employed a depth-first extension, one bit at a time,
    pruning only when the windows used so far plus the number of unused non-palindromic
    reversal pairs cannot reach the best length found. Complementation symmetry was used (first
    bit 0); nothing else was pruned, so every sequence of maximal length is enumerated.
\item \emph{Integer programming}.  In this approach windows ($n$-tuples) were treated as edges
    $\mathbf{w}~:~(w_0,w_1,\dots,w_{n-2}) \rightarrow (w_1,w_2,\dots,w_{n-1})$ of the de Bruijn
    digraph with vertices all $(n-1)$-tuples, and where an aperiodic orientable sequence is a
    trail that includes no palindromic edges and involves at most one edge from each pair
    $\{\mathbf{w},\mathbf{w}^R\}$, where the superscripted $R$ denotes the reverse $n$-tuple.
    The integer programming setup involves binary variables $x_{\mathbf{w}}$, one per
    non-palindromic $n$-tuple $\mathbf{w}$, with $x_\mathbf{w}+x_{\mathbf{w}^R}\leq1$; for
    every vertex $\mathbf{u}$, the out-degree minus the in-degree equals
    $\sigma_\mathbf{u}-\tau_\mathbf{u}$ with binary $\sigma_\mathbf{u}$, $\tau_\mathbf{u}$ and
    $\sum\sigma_\mathbf{u}\leq1$ and $\sum\tau_\mathbf{u}\leq1$. The integer programming goal
    is to maximise $\sum x_\mathbf{w}$. This omits the requirement that the chosen edges be
    connected, so its optimum value is an upper bound on $\ell-n + 1$ (as is also any bound
    proved by the solver). When the optimal edge set is connected it is an Eulerian trail,
    i.e.\ an aperiodic orientable sequence attaining the bound, which proves
    maximality\footnote{In all the cases examined the optimum of the relaxed integer program
    was attained by a single trail, so the maximum length of an aperiodic orientable sequence
    of span $n$ equals $n-1$ plus the largest number of edges in a near-balanced (with one
    start and one end vertex) subgraph of the de Bruijn graph that uses at most one edge of
    each reversal pair.}. The problem was solved using the HiGHS
    solver\footnote{\url{https://highs.dev/}}, seeded with a heuristic solution.  For $n=14$
    the solver could not close the final gap of one edge directly, and so the problem was split
    into cases.  Since the complement and the reverse of an aperiodic orientable sequence are
    also aperiodic orientable sequences, it can be assumed that the start vertex $\mathbf{a}$
    and end vertex $\mathbf{b}$ of the trail satisfy
    $w_h(\mathbf{a})\leq\mathrm{wt}(\mathbf{b})$ and
    $w_h(\mathbf{a})+\mathrm{wt}(\mathbf{b})\leq 13$, where $w_h$ denotes Hamming weight.
    Separate integer programs, each requiring at least 7752 edges, were then solved for the
    case of a closed trail, for each value of $w_h(\mathbf{a})\in\{1,2,\dots,6\}$, and for
    $\mathbf{a}=0^{13}$ combined with each value of $w_h(\mathbf{b})\in\{0,1,\dots,13\}$. All
    21 were shown to be infeasible, establishing that 7764 is the maximum length.  As a check,
    the same model was used to reproduce the exhaustive search values for $4\leq n\leq 7$.  As
    a further check, the same searches were carried out using the SCIP
    solver\footnote{\url{https://scipopt.org/}}, with the same results, except that for $n=14$
    SCIP confirmed 20 of the 21 cases; the remaining case ($\mathbf{a}=0^{13}$,
    $\mathbf{b}=1^{13}$) could not be completed in the time available.
\end{itemize}

The results of the computer searches performed by Claude are summarised in
Table~\ref{table:new_search_results}.  The following observations can be made regarding the entries
in the table.

\begin{itemize}
\item The searches confirm the correctness of the claim made in \cite{Burns92,Burns93} that the
    values in Table~\ref{table:burnstab} are indeed the lengths of the longest possible
    sequences for $n\leq7$.
\item The searches also revealed the number of distinct maximal length sequences for $4\leq
    n\leq 7$, namely 2 ($n=4$), 2 ($n=5$), 22 ($n=6$), and 3294 ($n=7$), where the
    complement and reverse of a sequence are regarded as distinct. If two sequences that
    are equal after complementation or reversal are deemed equivalent, then the number of
    equivalence classes is 1 ($n=4$), 1 ($n=5$), 6 ($n=6$), and 829 ($n=7$)
\item The Table~\ref{table:burnstab} value is maximal for $n=8$ ($108$) and $n=9$ ($210$).
\item For every $n$ from $10$ to $16$ a longer sequence exists than the one in
    Table~\ref{table:burnstab}.  In every case the new lengths are largest possible.
\item The status field indicates that the length listed has been found to be longest
    possible (\emph{maximal}), and how this was established, i.e.\ via exhaustive search or
    integer programming (as outlined above).
\item Every sequence listed in Section~\ref{subsection:aperiodic_examples} was checked again to
    verify that it satisfies the orientable property.
\end{itemize}

\begin{table}[htb]
\caption{Search results}  \label{table:new_search_results}
\begin{center}
\begin{tabular}{rrrl}
\hline
$n$ & $\ell$ in Table~\ref{table:burnstab} & max. length found & status\\
\hline
4 & 8 & 8 & maximal (exhaustive and ILP)\\
5 & 14 & 14 & maximal (exhaustive and ILP)\\
6 & 26 & 26 & maximal (exhaustive and ILP)\\
7 & 48 & 48 & maximal (exhaustive and ILP)\\
8 & 108 & 108 & maximal (ILP)\\
9 & 210 & 210 & maximal (ILP)\\
10 & 440 & 444 & maximal (ILP)\\
11 & 872 & 882 & maximal (ILP)\\
12 & 1860 & 1878 & maximal (ILP)\\
13 & 3710 & 3752 & maximal (ILP)\\
14 & 7400 & 7764 & maximal (ILP)\\
15 & 15467 & 15534 & maximal (ILP)\\
16 & 31766 & 31810 & maximal (ILP)\\
\hline
\end{tabular}
\end{center}
\end{table}

\subsection{Periodic sequences}  \label{subsection:periodic_searches}

The integer programming approach of Section~\ref{subsection:aperiodic_searches} works equally well
for periodic orientable sequences.  A periodic orientable sequence of span $n$ and period $p$
corresponds to a closed walk of $p$ edges in the de Bruijn digraph which involves no palindromic
edge and includes at most one edge of each pair $\{\mathbf{w},\mathbf{w}^R\}$.  The integer
programming model is therefore the same as before except that every vertex is balanced, i.e.\
$\sigma_\mathbf{u}=\tau_\mathbf{u}=0$ for every $\mathbf{u}$, and the goal is again to maximise
$\sum x_\mathbf{w}$.  As in the aperiodic case, this omits the requirement that the chosen edges be
connected, so its optimum value is an upper bound on the period; whenever the optimal edge set is
connected it is an Eulerian circuit, i.e.\ a periodic orientable sequence attaining the bound. In
every case examined the optimum was achieved by a connected edge set.

The results are summarised in Table~\ref{table:periodic}, which also repeats the periods of the
Gabri\'{c}-Sawada sequences from Table~\ref{table:Gabric} and the maximum aperiodic lengths from
Table~\ref{table:new_search_results}.  The following observations can be made.

\begin{itemize}
\item The maximum period of a periodic orientable sequence is now known for every $n$ with
    $5\leq n\leq 16$.  For $n=4$ the model shows that no periodic orientable sequence exists.
\item The Gabri\'{c}-Sawada construction achieves the maximum period for $5\leq n\leq 8$, but
    not for $9\leq n\leq 16$.  Nevertheless it comes close: for $n=16$ its period is about
    $99.1\%$ of the maximum.
\item The longest aperiodic orientable sequences are strictly longer than those derived from
    periodic sequences of maximal period: for every $n$ with $5\leq n\leq 16$, the maximum
    aperiodic length exceeds $p+n-1$ (where $p$ is the maximum period) by between 4 and 19.
\item All the values were established using the HiGHS solver.  As a check, the computations
    were repeated using the SCIP solver, which confirmed the maximum period for $5\leq n\leq
    15$; for $n=16$ SCIP found a sequence of period 31776 and proved an upper bound of 31777
    within the time available.
\item Every sequence listed in Section~\ref{subsection:periodic_examples} was checked again to
    verify that it is a periodic orientable sequence with the stated period.
\end{itemize}

\begin{table}[htb]
\caption{Maximum periods of periodic orientable sequences (GS: Gabri\'{c}-Sawada period; max.\
$\ell$: maximum aperiodic length from Table~\ref{table:new_search_results})} \label{table:periodic}
\begin{center}
\begin{tabular}{rrrrrl}
\hline
$n$ & GS $p$ & max.\ $p$ & $p+n-1$ & max.\ $\ell$ & optimality proved by\\
\hline
5 & 6 & 6 & 10 & 14 & HiGHS, SCIP\\
6 & 16 & 16 & 21 & 26 & HiGHS, SCIP\\
7 & 36 & 36 & 42 & 48 & HiGHS, SCIP\\
8 & 92 & 92 & 99 & 108 & HiGHS, SCIP\\
9 & 174 & 194 & 202 & 210 & HiGHS, SCIP\\
10 & 416 & 426 & 435 & 444 & HiGHS, SCIP\\
11 & 844 & 860 & 870 & 882 & HiGHS, SCIP\\
12 & 1844 & 1856 & 1867 & 1878 & HiGHS, SCIP\\
13 & 3700 & 3728 & 3740 & 3752 & HiGHS, SCIP\\
14 & 7694 & 7738 & 7751 & 7764 & HiGHS, SCIP\\
15 & 15394 & 15504 & 15518 & 15534 & HiGHS, SCIP\\
16 & 31483 & 31776 & 31791 & 31810 & HiGHS\\
\hline
\end{tabular}
\end{center}
\end{table}

\section{Examples}  \label{section:examples}

For compactness, each sequence is given in hexadecimal, read left to right, where each hex digit
stands for four consecutive bits (most significant bit first).  The final digit is padded with zero
bits on the right.  The sequences are also given in binary for small values of $n$.

\subsection{Aperiodic sequences}  \label{subsection:aperiodic_examples}

An example of an aperiodic binary orientable sequence is given below for every $n$
satisfying $4\leq n\leq 16$ of length equal to the value in Table~\ref{table:new_search_results}.

\paragraph{$n=4$, $\ell=8$ (unique if complementing and reversing are taken into account).}
Binary: \hexseq{00010111}\\
Hex: \hexseq{17}

\paragraph{$n=5$, $\ell=14$ (unique if complementing and reversing are taken into account).}
Binary: \hexseq{00001101001111}\\
Hex: \hexseq{0D3C}

\paragraph{$n=6$, $\ell=26$.}
Binary: \hexseq{00000100110111000101011111}\\
Hex: \hexseq{04DC57C}

\paragraph{$n=7$, $\ell=48$.}
Binary: \hexseq{000000100110000101011001110100011110110100111111}\\
Hex: \hexseq{026159D1ED3F}

\paragraph{$n=8$, $\ell=108$.}
Binary: \hexseq{111111101110000101001100010000011011101010100011110100111000110101111100100101011011001000101110011010000000}\\
Hex: \hexseq{FEE14C41BAA3D38D7C95B22E680}

\paragraph{$n=9$, $\ell=210$.}
\hexseq{00CE32BCBB6647B05BE6B9E04FD22E1269B8B10A77A1AA5A8C3FC}

\paragraph{$n=10$, $\ell=444$.}
\hexseq{0057BB6AEE7F7835BF2B9A096E45B1C27B3D7452C9C5F32D7051F62C3BEB31563287C94E9223A19CBB08D3C4DC7A6C042A6A20C818955FF}

\paragraph{$n=11$, $\ell=882$.}
\hexseq{002F5CB7F20BBF0EDF44B291EA7BB3FD0AC943AB754DCA8B1F9DF281336F4E36D27A5AF3EB0AB598B511B8BC2E05B2FD9A7EA0D35DC8A63BD91CF1D201CCAD8461E5CE4C4F347D8CED10A3C1968667B5C41C574C92F82C87CC08986F18AAC52C486093A36752EC04270D09581EFFC}

\paragraph{$n=12$, $\ell=1878$.}
\hexseq{0017956F52CF5C7FD6E5BE2DE4C713F62FC6F71C8DDF83A7D9FABBEDD4FF41591F0DF5EF0AB1E8F7712EF494BC08BEAA9919BD158B5D9ABD0AE7F70F5491E725C6DE3EB3BE6A737053E7B270A6EDAD35D518DC9F5027C67833B203E617C272FED6673A279B2ED47E5A2E9D2A8A77996D58704C5C334E2A12B80C665235E234C36D12CAFE20D72A5381A5D44D8AF2177683D82E496E30858A63476050DD1C23AB58D80B50C1B4A3C3164AC01058398A0914832B675D03B9E9EFA687F2500E3CA6AB904451D0ED9506225443664232E64D616C913BAFDF043C90EB6202AB48688E7D2064428B38D26134DFFC}

\paragraph{$n=13$, $\ell=3752$.}
\hexseq{FFF3ECF1F561BFA7ECCB7CBBD3BDBF96EBF6C3D7D2DEC7331E2AFC74AE631D1F479F929DD81F8C80F3A29D0DE95723A16E3C179784B3FA9FF6815E3943DBAB4154D4B6A3F7C6E2A13FDCFE2DE31384F7BB863E7A9B79E689BBE6B959EE6ED129BEB7E1AF36E9BDFA4F1A7C3B5A599EC5EFA9533CB97540F99DBE8B2D7794A1896AA8DD65C5E2CC073141FFA33752766B7A8E7719F8926DA76900BF517AE97987DC8BF38EFCDC24EBE17C8A666D374A659D950761D4E0F432FC1A5CE64FBCE0B6F8874E79381EC175192C30796A7AE177306E7A227A674C9567136AF7033EB323E49985669ED53C8D49176A63925EA9872FBA877DA85EDB075BBA47947BB712BD64E93570169C42D643FD68DB90DD01269518DB3793789FAC5C83868788B5043E3832D2906CA30CD851D5435B3AD63356E7579236873930D8D0A2D272C89304CAEB02C15C6CF9053F0A90B2B04237EE8E421B5D181A3BC3FBC1214D366869D44C47F990CE5565808F5B700D5780430F207871AE221A80638B430297FB208EEB9F4D251C281341989A8A58C2511A62D97A07484D8A56C7A454EC4C233C22C764A1562B1FD0980289194FA8456BD9C3A20FB1974207120451A4462C991092D812B54584882E081735883103280C148D1A097C009CD000}

\paragraph{$n=14$, $\ell=7764$.}
\hexseq{0006BBD3E0D5977559BD37EBEA97FF3B5FDB96F7D8DABEFF09ED975FEE17377FD679D5DFD3D8FE0E9E73A3733FE856FD3B393D7F3743F97CA7568F4397AF8F762B75A99E2CE7F6B5CBC79EB46CFC6BDF8CAF78DB24BAB0CC5611FD90A7EFE575DDC58FB326FB4C2DDF7CFDB19CD98FA38FF57505DF9AD9F73E3A7496FC3D9299DFC3438EFD8B352DAF15ED1AEAACA7AC1BEBB7DCBB23AA9BE66DF526BA53C64D273DCE5D4943552D5B7C5CB3DEE9F8318BFCD232FBAF57CF6AECD911B41997A3378BF779825042F2E916EA9DAB1DE58D73F5676DEE698EB7D2E8C97AAE6DB503FA9AF45CC7B0CBCFE4EFC10FA1BFE271BFBC46B50C5C07D6CB76D8BC6A2920AD3EB95F22EB9DD2DE5590F7478C3248AF23DA1533631DAC9B55F0E7E54D4FB79BE8D4E49AB8F5ECA64C2ADAC73B07BCA3E512DC47E4CF89A918A9E1D12BC4B587A6B8BDC4FF68BF32DBC70556BF2DA3A590969BC9A83C8DE0B558D8820E7A9E9E266BC38984EF5C760539A8FD9DCADECF9E4FBD2B70FF92D493F60097F1F16E0E66AD69A6EB8852F8AB14D684AE70ACC6B3AF35414DB499D371B8A667B38C860597C14B45998E3981EDC0C4531BD50A19E568AD9836765715DE44D3FF51AC2EB774469DDF2E7A3FCA4B32C963CDC17CE80EFAC8ACB656D08DD6715ABDDD9451EE528BC958C4F9CAA8D8626DDEC52E4B8D962EDAEE4A441EBAD2FC49F455FD137BD0CDCF112FA13035BD23B34F01869AC17B6117D45D8F0A819AC83E16783EC463D5C3932116DC7EA1F9E1BE435C1A148FBF555209FD05A493AE01B726238D32364F9616A23EA02D1384E2EA8CB544F28379DC67A58C08EE1BB989143B67A09B922589DC19DA7AA8AE830A694D152ED295DB12E68396C70D9A1094BC424F7D012B60E407A4AA7CC648C15E17E11CFB5543AFC055117838B4EA714AC1396BE9CA62D87DE18DDAA686D6E81AA4E05EA543720A101BA4D95128F260A0776F5BC29542759B4B0EB40F38F2163A8C3B274415142E99CB2B3EB13412EF8E480CFA4F0965C481957812C2091103E6889FA6CF20EF009E90D94817146394670D4C3849BB0499888A2E30B85B3846AF322C81D502DF0A00EC490D1F86418044ED047C018E022C370892A62057330952CE002E8B1508AE4E22D6D82C405308CFE2813141615D84DA29D8D1268E723233406A01C121F20454781D215C118447476E98C2468D0C3E30402932E54A60F1DD4CB43CC59661E200B038832E19D6C5131A0417287021D4B84CE255C8E5813A71D0B27225477C0DA41E10DE890A3842B129C306442E452C18B04521DD092C072A0D3002A4C06C8081B0A34AFD247154C82210CE420D8F928C0510F403450741FFF}

\paragraph{$n=15$, $\ell=15534$.}
\hexseq{FFFDFBA59DDA3A75830FA6BC85166E9D18F43F8F6FBF9F33DFF3AB2F57478FD9DA29A65C15F75B14C53A67FF6ABAE7DD9B7F2B2D8FAFDBE2FD54EC09FB778E42F0F8BBB6B4BF8BDFE8E7F4FFD97FF92778CF7757EE6BDF3F7783B5FA9FA973F6B67BC37C3727AFC8DFDCA87E2B9DEA7BF36F4079DBFE955B1DACC4E75F0FCCEB2C4F1F5781DCBF12ED8F3FF0F1649E416FDEA0769F23A61BDBCEA71F69A39907A8B62775AFF6DB2893F095FCD6D1A3BD3B7BD36CEF03740B2BD6B616FA1EFB29517C777A21CDDCB7F0EF98A52DEE60D067FBE3F009D266B92E2BDAA7B1F08F2470F756F9A6ABE7C3AF7067B643B8ADC5F66B47FEA6EBE8CF58ACCCEFAF65BBE8B2EAD7DEE5D4726DD99E9A0EFE377E9E27FB4674D3F7272F49BEAE65EDDA55F05581F7B5A92E9BCCB68727E984CFB3C8FEC9FE72075CDEC685B6BA7F5E9EBACDB72BED7D47B8F2CF600A6AD97230EF66E2D09F64776126ED0A37EA34CE17BCFB75E276E9B035DE7AB8FF37FBB9332EFB5176F06B3D1D1BDD72D7CBD9F7EB3EA66F527F76C1F86DCD16B73EB34B7395E16E523A2DF2D6B0F3760F19952B57F96CD83B3263FEC7BF82F722AF3797D0A07F335971F30343BCD45913FA1DBEB8395AE6A774D44279F98EBFC6BFBD82FBF25E7F140BEAC67CBBA94E8AAF287B637E64FBEA8D5AA3FDC6F46DFB0D43EE719DD44FCECD49A83DEBCF95DF1A95CD3BE4589AFB0B8A73EC65F274467F60E6DD4995A8DE79B69EA610FDB54BAFD2F361BB704F8F5F3329FC33B1557187A5749EEC07B5E976BB6D2F74E2E86978D9FCB1D9F94A46F183A54F5A8A9F8E693B9BAB1321EE047F83167FC48CA58DAC65764ED96B22F4304F46BB307D2CBF43631F6669F4A79AEEC8DCB9A158A0224EA9D7AC9F0696A26316A7D823E303CAE1DDD21E5DA47A2B9655DC6BA9B7C96C231CBCEE8B0FF28FD721FF595C897FA8A4D9E28973B69B3B283C5689D0E9FA315E74FBA3659A1E3ADE950D7E1CB97C5644371FAA67C0F136C73E9F115B701A3719BFA848D9D548FBC8E559D859F21E1B2B747743EA55D48D73CFA27B2689BC5C7B38D188A7BA151B60AE6374AC39FAD1B4A375C27065E8CBCAC6E7BE90A919E7098F34AF48E6C06797932943B1E18593DA64D554A839F1A6F83A0F4EAEDDC117FC86FF0A4F9DFFA0EB4CC3CA6D7C4E5BACA305A97DB0AAE8079145ED0967E36226E7719426E9231B8A835EDB0C7F279985F10D746A6949A0BF78060AD3C1F1357E6AB50EB38ABB1D4671C2EA1E896DED4F2E324E37223DF293D5B99CFC147393C6B500ECB2D484C4BE6D21B662EC49CF73CB7AEFC6703F714BF3DA06F5B0138917E16B11F62DC36AF9578A968F3A3F3FD39E4B0DA19697543CF2AAD31C7F5B79017F6F549E1703876C4B45E78945817AA097592EF1336E2847E424E06D3A135C6AC92F02E9EF7E18A3EE258E2F512CB2750BD92487CDD12AE2E421ECACC3600758C678F4CA381999175E8B9D6E483F6062A28EBCD10A72C5FC22AA6A863CA82E4A41E612F614D8E97D621EB9815706A49D3D22DEC86CCC6CD39CF1C4EB4478D2E5867A0C954277110B674BA30097EC82562497336B94D6F80C2B2B4890972BA2805F09A962CD4225799283B0C3A6E5A6EF6C08A668FA889613B4C4D9BF7089158E5069042FF0327334265FD770AA3CA35BCDCCAC9C866F044CEE1AF0D323826D301CF6692E01F58073568E46B0D1C4ADFA297701371D23435BAAB1ED2640C6D07CE830C6A06C590C5E3B6F14DEF16F0CB39C3C871B04B7E2C66CE5519A46292BE99B2931A76A51F24B71CE0F02A6FCA1739862134F6E3143C10CE8646C38F81365F7004A9D023BCB4842D5C82BC688A9A32B0BCE6FB4B87AAD02B8983D007ADD8BE4016BC65CE2698C1BC49B461AE2CE1D125A539941582EFE9214354A161D5B4E3A17E61CC2C2D429A1C34A12AD674066D111AAE50CFA1297A0A5E464BF0716DD093170A2B0982476E62623B41B22A4B18228D261F22128D150E7545AA6E270DE49154B262B62003403396ADE521DAA237D382EBB561046EEF8276725AA169AD712571C07EA03A386558B4FC9851E90752DAF38F18E8846C84E620F86AF62190ADE21247CF0613DCA2E08577AD2538A4470AF921D207B184A8CB1AA31A07D1EB30242C5CDC3C009AC7060806E210521ACA93B25406FB902F240E5A0384B1227A865392B6C968633575AA8EA65A2056D70CD8AEC14530CC64F058FD03A73411C023C222905E2C1C318E21D57BB504627D7505B53A9035185D97848B4644A9607510192533C4611D74C93635E8473C3D840D3E163A09A119B9406250C4816180391021F6D020E678419C46E061CA7A1658570B41C7AB3C0889E504CAACD62C4660AA63E9388BCB116532A38EFBD2508742285204F17808C0B883810869E4D85134906AE23D4A07198B0A5C9887158832C1B9C0B45199C61016CB0162A70A6453808582CDAA01331261255906E1916D8C83340CBB0E80198441A25025AC1090E58272B0105A5B2A906224A567C5B06522A053E821C8207E5E20986120E4C91663E54C441768045D3103222EF7180059E6460CC9203E5156659854106D4C14A2185049AD9339B50DD685D220914116116EA00409289C647E31340524B00A5AE8430288E94D48A599080482C42AFE028705081418E1C947890D05190234B20556193708032A15EAD2B6449A7101300DB5814C0C04B259480ACE45CC90000}

\paragraph{$n=16$, $\ell=31810$.}
\hexseq{FFFE256EF4E2783F628BE6FFBEF53FEF5D01BDF8EF5EAFF7CFD655CF11475CA4F77B530775DBFB27F7F84FFA7BD8E5F508AF6FCD28A3FDDD15C4797DD4E2FFCEF8DB752C7CFFFAFD5F7DDC57DD22BEB5BABFBF3FBDBED14EB733C7ED2F8FFEC8A9993F6CDD5D4ADE62D76B877B0B3DEB7AB0EFFAEF4AFB525D5527F43A74BE9FDBEBE1FEB5340B8BCFFB79789AA9D7A9CFF07C9D05578B665777AE6DAE7F4EEA4ABD6052BBDE891B67FBCFE2B7126B7DF8BFB685F32EEEDEC44B9077DBBF7BB9FBAF6DFAF7B6C2C4EB62BD1DF54BFC6F4B67DD6DED7536FF759C5D7F3CB7FA4D239FB55C7AFBC11D5396F9BB547E8EB2DCC5FCDF4B3D8C9E7D1BFF28CD974F8DDF5DD3BF4E1F92FCDABDFADF813BF9591BF91EF0160D640BEFED73387BE173E3F2FEC6BC7F7C7C47232DFBEB74DB9CF1C773B3F7A4B2CCABFE91FBB3D75C4A735FEC29D5C18FF871BFCE9BFC138C5DFFD5BFFC968DC632FFE6EBF77D0AFD2EFC15FF31CFBCDCA72283EFD4CF9DBE5EFA9AECCFDC55E80FF46DB3D4BB2E3EEF0FBE69116BDD5ECE376DB56EEC6F9E337B67799C877E36AD68DA9DEDCDF86E48FCA6C0568E5C6DE40EE66DE3B88E25EAA53E7BBB9BDF6D67ADC6C62C3CBD17E5DE23FEB19F7375357F1D7B338A3DBE0F7B90F36FAB6E0C07B755CAFCE716F5F562B3ECF6861ED5A734B754C166E39B76BD45EBCFD3E879CA6A4EF51F907B6EF9DC97DF17DCE2302BBF856ED1AFD8FA7CEAF10E57D56D7FE35CC55A9FBE46F55F0CF7EC802ED8F163F3B59D9D86ABAFC7F33A76F1BEA09D4FC2FD57461F370FED8CE7BE9A13B760F76CA0E588EC282FB714E73AE1B3CEBC6E9D4D386F5AABECDAE2AA14275F6132F0786FEE9C7CA793535BB1339BCE3F0D5349F65F45AB33EB3C6CEE70888EBCB76EB7D3AEA33C8FF0D0EAFCB7855FB83EBB6A386BE54E0EBF468889F3997F5F4C5B2FDAC4FCE89706C5BE4C4C4F970F176F597CF9B71D3FC59925BCEDAFA5129BDEFA2B5E17B4EFBCB605BACD3FF84409CD41C0BE44EADEC154B3B8A587D95FE613D57087F5E1B7BAC76D3C5FB3C1DAEE3E8F7FEA5BEED6FF2DACB86B05AEB30E5589CC704BD3B5FCB5BE3C9FF15368739BB2C5ED8B3D3E1D7CD8C8F6D558E76AF3EE07EB67951DADBBAE45876DF48D5CE4EE0E34337F5C60B4BCEE129367E725B7CC3B4FBECE5EA07FDDA6C27ECCA1ACEAD7C6B746BB64AF2EE79F973D83B3AD0B866762F5B2B4B6E8BE0BAC4DDDF24EDB4D4E669EED17BBFA5ABC0EBD3F68F38D944E9F744DC68DE7AD670EAD3370BD4C4B74A3E7A7201F7B8D5E4DDBCD745A3FAC2E348ABBF453708ED2EC5775A619F939B17754AE7E8C675418C7CC1A6FBA1FF41A94DDF0C1CEC8FE48E4BC7B5B3B5381D2B70F659FDD81E5DB7357C4B6CC7772D72E5F832C37E92A6B3FE5E14F5F99459DA3BD379B45CED27B70B7F53E9336C20FE7666753A307FEE7BAE93E576A2DCAF3CD638FB02A2D8EF95BF9C747FE463766F4F5B6F3F16EA3A997EB4F4F3B43F975F27E86AEF8AE9EE520D25A57057347B0FA73A3D7B18E372F9AF7C6266F1575AFC94BE430EC86DC217F50D3F910F494B5B463F84AFE42EF22F5716D3B4EA93ACA67BC7D7A5924BBDC0BBB4A9FAE5A9EDA663EDD329DD0D15784F1D6F0750EB26FB627BF2435EE687B1E469EC4BB575D92F654C95D7A175D6C9F4C9C43F6EDA27FB9544E2E6BCEC7E0E9E971E34EF1E7B93F9B028F1EAC9D73524D9194F498F7D8C24C83E7993D0C74687969347932BC546FB317A556C4AA7F1B14F3DAF9EC93F2ADE5A126FCF14FBD0DACCB9CC2FB29D4630D8F4E0CCBEC487E6F473E78456F452E7928BB34E9B5A2DE53A78D62FAEB0DA94A76A01BFE3935FDCCAD96D5C77ED857712D74DCF0D959B6F8B5531BAB61241FAAA83BD18AD07FCE156C91BA0543BE201B379737C6840FEF2BE72ED08B260F284D7AEC76064EF90A49CA29248FBC92A29EC324F5E01D583980537F1A6DB88993E32752E50336E939137B9D381372E91EE7C2AB3F87ACA39D60C958D757D2447DA4C3E3D10A3EA617A4F5367BB165EC25F4388567B8A366D0AEAA13D946F741CFC32F80CA71CE5B634F998E9912C39D450B7C39F41D6D2E6D3E2F7F41770A73074563A474AB7259E63E3B7A67949A0D1C7BAAD632846F1CAC1B3A93EAA240CE6652764ED07B98D143B587C16518F5A07E7C9074D8D6D88F5DC9A341E7FA6EE5751657B48A7E9ABB61871E27B4CA87E90B83B1E1FB925FD13B1827BD31FE6856CD78A1EE396637E63357A61366E64DE04CFF2A22EE6B8E851D979C5650BB4E32C8FB251E9B63912B51A07BF0B245F874057E18D93C2BA12569BD63D53ED7894BBA1637B2EF81ED06FE0E4FBC4286AAB9BC263F65BDF33D551B2CA8BF8CE056DA52E4EC55C9EF3A853247A139F13EB4B8BBD5E69E6AE5ED6A9DB34DA745E2F1A26B66257D1333DB2D791DCB8A8096B5DE12F3793E6B44E7790EA42AE37C536AC33CE6D955FE20E69CF6F403B7E84DAC1F98A5ED0293E2B132CBF890D367359BBC63ABA5EF54D88D18675FE4C3AC559E8CB506D4556A26A335F5C58559910CE1431AE35A6B5A8219B4C812FE52F3363AEFB8F37497C46B43EB117996E3C35C816CE5548A22777E5945794A3C12DFC523E475AC59B3E3A8E5ACD59133EC3A32B43A1E31AA5B4373AACE9A8B154CE5C9B7AA16AEECE889A338F4C498D397FA8ACEFC8E383434CB959C04F1545AE125D2D39292B0D1F356FB01D1DD2CC866AAEE64578E035D227E03D4F391DB947E3D9F2A07AD87F4A6DAA2A74874C75A56DF8A0E83D70B88AE4B969BAD59694F0F8A331ED28D3A81603F90535A0BD8B9E0DA5616ECAB68A03EC5B7A45A9D621E2B9829F2D57A923CF21A5D78437E5258F70CB39FCA03195CE9743B80C57F053B0905FBF01365A3272429B265237395278CACADD06EE9302E6CF170AA7919A9EB25D15652DC8E704F0946EC2F853473196359486C179B13802C79DE909E5C4E962B852FDE3A2E81BA91D983EE5132AADD7900A4B9A4AF9766353B2B869E8B34CFECE0F5865D73EBF09BEB238AF06A1724CF66039C837D190C7E983606B886F302D3B351EEE18AC3939702AFA12ACE776E33DD2A9BF560281E70CAC7084F989DE19BB845FF6895284BA9807188F191B4E12DC172AF9A3A69E0ADF0243B4666D5686C4DB47C1A53A21A16E05314C7D51895185AFB4D35F17513ED916E894A2D452A04C0EA78F093BB4161C20BB2864DD1C6C13149D345DE06B18A1C654BEE1F46010FC4E592AD80D4CC54C45D704DB0269B91B8BE833D4655A5555A4340D5AACDFB51933F1D00D2D8B8236EE3B1D125939E8AAA7C4584C15393664186FB98789986EF66B965886D3417099E2772CF830B9A0E45D3CC007D5C9265D5A97E452803F40C77957566CA4C624B8E28B41DC6F2460252CDF2CEB3F548F033FAA2037642BCE7DAA0AD8DE241EAB1B4C7A952EB77040BCBE583AC2B5C9C6BEE456152F211B419A62AB90E14C368E20AEDC12BC3A829CAE390AD5A047374223DA843670DB5378D4BC3CF401B1B7224E68CA5D41AFC589698DAF819798DC4DF3A4DA50B4A8F10D508C3D5F299A91A676E87371643C17E609D9C88F0A8FA6CA95511FB0662E65C049E919E1DD5CCA0DC56B31C950E0F8063E847F7194818B78BA567865323A9485BDC4C3D35928543786A7841B2F5A4E8EA12D9A11D1C143CCD152B1795D144F506614B2F2ABB0F7A0FCB162D8A56B63812A78A244F4CD85D1DB152985030AF78DA48E823E1607CAC53A1C3B07219070F40FA026A581DB21D8D9ADEF14CC258DCD354F01C37692E182E9CA5C62A2BB223429F9A990B63007CC96AC7838A93736452675DA60A0A55262E8821F5288AD3CADAD837A91107E144C58D8929A6046B089371044909C2F3A00399885CC781700904FB1051D102196D20785B60E190E47A87DB06135ABF1AA39B3049B1610931213107D40966A99468EE1CA7FA0515496D44AE54D188DFBA02270C09FA8F41C7C10131D9E488AB052963443E30F248D7A4836C50DA0393095921CF5C916C9B82922A393D62929831CD58D060D564B1987F06F9489B415EC9C2059CC5099D2A25181B7847EA4E6480225B0EE50A3493E1B0A4E090831E941ECD13C4B175A9350DD608E52BDD92448F8B9793CCC10C34EC12C8BC95B42153D5D121C733E218C08B71875452144BF447A443AEC1A46C84A801C405668E02974CCD628119B9338D326636114B51E630B0B4D2E0B0A28B20F32A47037A22CFC467F3EAD1A55C08C290C46611C6602FE0165DC0C4F47F1C44A48655389418B3160F103E142FC70B47680BC76542F30B608D43D402FAC02686009F269F055905EB33A0C914C407C8B9325C012EE426AE246A6EA867CD01A7426176083D2C1719268B6D9311D27963BA3F036BB11C3C88C89E7109BD0A4640C889199019F51AD88B681D6A8142477FA2C322E1C5F28619B0517309130AE47D12D4AA38C35061446E5432E626B3064152CAC272B5812E71E5817BA602A6ABFAD1E0361BA6BA80F34CA02C9EE0D7991CD21F7108BB1A8D87180D208F3458E546A231A6B0C05A908FAE7759271C64594A1A54901691266043DE5045A7D798820E524C01ADE4C118F308732D160C550D92D94D2D15F1CD12CA049D24D5B50CD2BD29579644955C5304D73D504A461CC49140E28676C160A8D30B414C0C6A124E1E020B6E590D6D743120AC5A0187852B2F4189A58E2481D0CCD26F688A0B3C8234A0938B18C430AAB2E8DD420EAAF70B35809C500E510BCA38BC0C19DD887169C0765B2042695A309D812699A09E820CB881BC883B0C538C21AEEE01BA4303804582215191C26A0BC19E0736A40BAB190A912F9D80CE918EE413382A4A1FE405192416A910B0D9A490F870C64103011EC8C5C5AC2CB26C1F4AE9A42CE13E1389B1F038ED562E437467D2F1601E42D6E404881CF8148904548323080E0C4B311488338915C07C4C94A9DCC4453C607F111B020344A1CBAC882F8B2D51307138F00A83D1421326B128B030C4E8A8865141379530C18890BDCBC0094226C339335C00D547C6A068D5B80B9109D314560916216A1C962611B318373035492D0F46C4092D64FC6681C82EA4046E1B8FF6E4F22A18BC41845B002DFA888432505D0AA417381A917869A230D009622A08256232EC835A5893819A893E8808779C1C114D584AAC44370ACCDF04A0212C2F19F0852E140554D42AD10914824C2CC349E084792129EEA20DE30AC209C9E04564D0BC242C518803DC0446272105210824A394205658021A48465821688EF1694081E901E87CA2E26892571A0A8480DBB32694CE22E06C825C5204142CAA0654146E0009D0DF08CE58464B0222FB413288E08748443144382DAF6AA045150B0011DD840EC01DFDE1754806422A8046825A2643E9050041AB1085D31A048D010032D04608BD08201EC56C8F29A1CA001BC18D2FC4C0610B228D501D71C390341142B494400A210351D080C98E61C4800AE33B298C9AA51C2BE40A0F1F20D09A085069C441542017483141806A00CD107803858C3201104D265AD0055D3A044C6664A90E041816468B8460000}

\subsection{Periodic sequences}  \label{subsection:periodic_examples}

For each $n$ with $5\leq n\leq 16$, one period of a periodic orientable sequence of maximum period
is given below, using the same hexadecimal convention as in
Section~\ref{subsection:aperiodic_examples}.

\paragraph{$n=5$, $p=6$.}
Binary: \hexseq{001101}\\
Hex: \hexseq{34}

\paragraph{$n=6$, $p=16$.}
Binary: \hexseq{0001011001010111}\\
Hex: \hexseq{1657}

\paragraph{$n=7$, $p=36$.}
Binary: \hexseq{000011101000110010110011110110101001}\\
Hex: \hexseq{0E8CB3DA9}

\paragraph{$n=8$, $p=92$.}
Binary: \hexseq{00000101001001111010110011100011011010101000111110111010001001100010000111011001010111001011}\\
Hex: \hexseq{0527ACE36A8FBA2621D95CB}

\paragraph{$n=9$, $p=194$.}
\hexseq{024D3A263BD423D9761E56956651AE17E6DF219C5B9E9182C}

\paragraph{$n=10$, $p=426$.}
\hexseq{0195DB57B14B78BF2FBB3937CEE8AAE6368C8EC1354E99FD470D96B30BCD38F293DC2B1D07C8B284BADC87AE0C49605882A4435F114}

\paragraph{$n=11$, $p=860$.}
\hexseq{00C3F4783473AF2FF05EC5BDCF487DA9BE84FC55DF9BCD0ED96BF609C6EF23674D46AA63BE61772BC75B19CBA94093B4FA8B49981D18416207C9581722F0AE1E7F7898A45386E652F8D3AB3822C92E285B614D321AECD6A0CBD6443896193C62423B203DBA599E54466D2AD}

\paragraph{$n=12$, $p=1856$.}
\hexseq{0056C5FE4E75DBCB8E6DAAFDF63EF3498DD2C986BA735BB857BB0F8F5F8EFE33B3F75A447FD66D7D3FA894CF58E2E7B3AAA230F623823DE8FC8F9E945ABD0E3CB51B8DBEAE32BE29C613EB8AECC29A6EDA74D729631575ED70CFB937CB327BA3208B12F9DE0E4226B16CBFB051595887F399B252A7E8560C66A4F80CBBBE83381B4AE806B412AB3E505C97B29E6292E6451A70917486F706CF165A3B9A2F8590F127876C8D85C42ED3313956B68609D80B72F0A026E05E57224B78D09950D5034C407507454938A10CA21A40F523B53BAB241E71307DA27AA5C20D1810547A22052154D8D5BF0587}

\paragraph{$n=13$, $p=3728$.}
\hexseq{003B9B7E6FF9F5DAFF72F99DF528AD14C5668EFDA64D5D8FFD3F5CFEA345E7E9BF1CEF79378B1BEED0FF6EA2CE5E506F08313759DB0A777A9A7BB2CF7F2C5AF25BD84CEA7DD3EAC667B3A6BDAEFC99BDD1E7C6A5BA46AB9C8346F4BB88BEF6F91983FA4E631C5F8AE39305EA27E716E5B8A18BDFACC7E6587DA227B525ABD95B3ECD6FE8A9CF7253F0F71B7A197662EBD3D48ABF423EF076B4F9AC16CBAD64AEF107EF9DD7211B66E57E12CD879DA8EBC33F863C54B8303F13CE1E3EC0EC55872D1FE68724877C1667380A5FB122D0DDD07579825F5A86B643C1A887D41C4E965EE09AE694EB6F3CDDB903E629F12DB871E6ADF4CD533369E649FE80DD10AC8EC332367A5708CB53820722E7AFB205E433B2943D31B9D26E3E1D697C62893CAC247AC75882D553D8E6E17C04362E4F1A13E15E3A42EDB572E819C2710D2D8C212750BC657586D1928480CD3A0994F0BAA58E018585CC47972AF5B14B4D98556551CA206CA3D24DBB86FE225E8960CA874E33C02E3024D3035F17828B8B86522693229239940F130C450B2DA804C76F03860AB8956C1164EC1CD4CB8815D290980254DB42E05A408446920F202C854603D54AD50A68C4A3850488648C88250589CB02A4FA35A105628160118D35218689B1ED}

\paragraph{$n=14$, $p=7738$.}
\hexseq{001A46F3B5EAD88C8CED3BED777CDAFA7ADBFBAFEC7F02F797713BBDEBB05FF071FEACFFB2FBE4C4BFDE6FE6EB3CEB75121F7FD6F63BBB5B2DD2FE3EF698E7CD6FF23D8CD9CFF47737D6417DB7730EC7BF3DE224D2AEF88D19F957BFA74D5DCBDC641FDCF0D1FF2C9BF05527705B78F8A9E0F92BFF9E311CB6B5DA177E3B152DBAE774162DDE470F6B4ADB6793774AA779944BDA3E5816E578D5BD68522A13EDEE4F177B2B4FC69E9B76A4DC63722EEA7F89CF758FB667CADF3FA26F83571EB30DE5BEABAFBD2BEECE9DBC1659BC653EE6AA52E5417AE6DCD8874FE28027B3E3A2AE99FC9141E344F3E8DDBECA3B58D03833DA5272F365C73B23183BFC2FD4F7159B55CBB6176399236FC8ADE4F9EBD8B5C3F7F295857F42B13C258DF323ADAB5E1F682993D51AFE437EB8C36355AC4A87CC07E95216D22E08DD5584D5A6749FA057D7E5FBEE9B962127BEC53F721C77EDC75840C2168BF17D21A73ABF778E36E24E3595CFD525E4A3EC3D3ECC4975DC374C01D79129ED913AD4EBE72ACA667A3DCE797AAC9FFA68BABDBA3587F3B941B7A39EA408D733583155A3A9C8457BB1B6E41D9B2678AF03BC333874656657CA95E8D33749A35242BBD1250F648771E696E1E5C9B62E53CC78B4998E876D533975A7597B38FAF1660D9A2A4CB2AE23CE1AB3494EC35F232873A182E96CE12B451E1853E1C8D826F536804C3F383A597474E3C2EA61389B0E34BEC1A6B6CA15376927D3865F99DEB18A170AA3D173EF08F0342756A3F36124ED5D4B05A729B84371652CCC5E16F14EA51CA2B84C575819538C1519CD7E1471C8A5218FA660B0CA4146D1C9D470B57D0D6899A0931FA8D96326BA31C0C7A6411984827C6CC9D9E087E25D3E443E33544D67543373206A9116FD187D5E6A5C6C583967EC4591F517C5C55BA0674712D34E840A9426E8A6BD27080F4B54139B88313A44907A0FA5381071630141968870C6D75C519A78983DEC2D83C2134DE68C0A2F4C53061A83ABB700BE163B458ADA325F34091CE49B532E7A51A516038889C9A1D725B1E5110DAC04E88FED214431292604BC4A30B38A8AD20D7988A2AA6364B280B9A25765AC94D90F562898CFAA06CAA23BE88AC039A87985EC0609132A203217442F1E419382555FC27A1CAC9254D48E92D7C05AD4319701E7A036A11D07308E33881F70A110432230C940F7484F2B1D438290837802C814F506A68460128A1B0047ED5542C4D45B607408BCAAF10574A085588807C4214E69886883E1B4FA12349E04E024B0B0A5A1961D5407761500DD9AE0DCF8916CC0FC9CC105C903B02269016089676264C282B411495C0ACB640412E18015CC6A31015B83D74}

\paragraph{$n=15$, $p=15504$.}
\hexseq{000CF9D866B7DFFADF84C94DEED4CE78925EF679FDBE9F6EEA6B6FDFED6CE3DFD6AAA6AFB2AD6E64F8A530FB5DB57E3BAFC4E7379CE07F71EA2EDDA5CD79DEC2EA74ACF8177B5A0DDBDC3FB4E535875DDC97A776677F87CBAFBF5F3F553FBE3F6F3CBE61C7DA6DB4AD3C475E5EBFA6A7EE5D81FB38BD746FD6DC74BBD3E47D874B1BD7DB3B68CABA975C7262BBCFA55FBB2D8CAD7B7C57F6449F6A40FCFE372DC53B63DEB7AE11577CB6F9AECEE3591F2BF6A1BBD463FD86D3D3A90BFC54FC1B2FD71A516A3F386E38A0317FDBA46BFD11BDFF374DFEA8F847F2FE597B27FFB86DFE576F1173A9EC59EBDFC75BB15F76E63664774E12E652F6E5937778FFCA5766FD3DB8D95EFB8C9ED2AB80CEF1497093AA7197FBADB743D04AF13CA7579CBEB63AAF736635EEB2F3A53C721F32615831B038F4AF9D654B3C7CDE16FFE774A90EFA6F52D1EFB335FAB8790F68B3821BA33EC99BD1E8B2CEBFC876EBDC55BFA5B53EFC983FD4D39BBCABCDFC9E313FEA0E74CFC675D5AEE9362FB46EBA17DF33CE41F8BC8BF0DE87C94A025475CF7126DA17BF75250FF1C67FE9DA6A1FA95716468FD0DAEAB395C79F2D6D9F77433FAEB220F9F7B6E4FC4AE9B20EE50EF2583DBDAB3DB52FF8673A0E36DE1E07D764756CB17706BF33B353BE734E9AD7DC8A2FA9C4229A95CDF41895AD2CC27758DC5BDE152F576B24F66C56D89F7034F7385FDD216CDC41EBE3D0AFD1F7A64D1AF1665E37FE12C916FD947F5941F7636AC2F5B4EBE865EA4D7198A3EA65F2736C6A58AE31BCDCF6B9DBF456AF8DA1DAF9A4E81A5BEC879DC38AC27BF429CFD8CC2D5B4572B7A7216B1ED8C49E84FDA437E2875B35A30DE3070E716D343AA36F8355F0B376C44DBF31076A3B281EFC6B9995D788319F10A6CF2758634D23B723C8ADE5273E354B0FE643105FC65BD9A7C66E0BF3A3978A24E1199FEF3841B986AF606BD23E0BB30F7C9ABB9ABCB3F94795EB4EF4B0B91D83C2F9639FE9238F0E5B31687F563A4FF3D845F2182AF1A63569BE348CF7D534BE3B451CF40EEC47E3C1F5146C2DAC5EF2EFA25CBAA19D47B55CB06F9E44D3CD627AAD027CE271BE017896EB928C2F36EF04C8CC1EAA2C5B148A5DA7FC68EF88A3115C12FE0EFC415EADF522D1D17BB866685A7E7AEF136BC09D9D3B32DDFCDAACBBB01F4E8DB015BAC1636465A84DDA203CDA5A4FACF4BF711A9BCC32E3E4478BE983B09CD04FE0B773A4CF77DB8399D192FD03394A58235C087F16BF056B3F7B923FCD518263EE8FFB312435D32C4EE2ED0C330164A40BD527EBC81AA16F5453EC4E519E686E5795AA087AE72696253E93C5D5F48C5E83D965AC85D8F3AAA54F4A3807E9C07A9890E021DD43CD44B6EE316EA26E53112CD744E9E4CACEA80E7D1DB0AA799566B2D29DE6C7E1C2C1CB33E49DAE910D6318767D02DBD4236B5888AD2A32CA66D4542BAC54DCD2E9A3E9A0CB6C30F3365A3574BD14E19081B9120E563CBCCE6C92DE888671DEEAE71187C73C16639DABB0F625A742C70F74C0DBD8F8B938A649CCD489FA6777E787EE854F0D9A536994FD1A1E5D1C6A1079629EF87AC9D3178CA93B143653D6E9FE05594A8838B0F4DD543EE7DA091ECC08D9D6A334D89BCB4C3A31A36A0BE125D95C998F7054D19D0E9B9361AD9544E62AB32A281D15692B64AADD00ACE362D5286A796F1D109E9885E9B47CF49A6BD9D8C2562E378B89F943197438E72D5408BC3B6CD18A8BCF7BE9724AFE4590EA1849AABA043554731D88786BEA1CA8072220DBAA472D17C129D16C26D999ED039B0917A8851D0A32856D72E30C9DFD88275EB37C83855B0FAC72E511A77313C18B053E02792B8A4346932D052FC601EE41950D38FCBC1D3991708419B791C29717061421F67E9143360BCB8763BC227BA2534476B02B4296C093CE0CF215468C8DD2D5CC82E558727AC2D2F094BA29374416A725D47282B8876009AB133C552E070D8F1EC96E7E60962032310DB28B34CF5AE251595B075609EC1478ECE1F4864B58579818E660BAD33891B9DFB0480E834DC2B52A0B229E0665DB278431F466512057863C4CB801A664944FABC29A159A01BF898DD41B4183E425A9B73CF02CFC90AE7F4709D40E25C484CCC652152669255B8C14ADCAA88E4AD1C4D2B8611C61640A59715E4A145CDA27130D8322C46439AE502E224CE18ACC63D68720AC988161C40F1A062217461142445B6A9F344B0DD9E21A4BEC299F54226610B9E5098C531D05161A948F050B0CAC0F4954CA00CA0894201E2DC850367E181F18E4915F973E1A8338504DC688A0E0485C1498F31E14530A424C9CB1ADE076949268AC604DA4E0F874302A46D20A874047D082C5041AD59613045811580E322A77410301E092D119C017347204BD8A835246A355309A9A0544A9CEC9C0C9B4023A7D40423066AE4600EF40690A511313875988F4847BCD308D09261E2042C03AF95322806189AA41A0EB4898411CD2064109541E8EB01074826990B3BEC0B454908918547AB44D9B7017480DE029E89A8F09086868912F2684770BC010C04C52D7C48AB109B00118E1104A201ACC81D920C0D0143417B65CC1A2EA004B407E171C1E61208018375A00931A190223837B202DC34C068F3D4A35028C0EA6052440E525946FAC83624C036B3B422D722027493C00B120148C937E58608125320102832564B21382D084E892C67A37009C9052019A21F9}

\paragraph{$n=16$, $p=31776$.}
\hexseq{0006EF2B25A52F68FB3E7B3D5F7FAC913F3CEE6946554677E3FFBD5E394FF9EED78424FDE959F278DB658F7726B29F2DBF719BFF15EEE8FFBE37F1F3FDBEEB30D76ECEF68B53731F77E90A16AF39ADDD8FDE2AAC7FE4791F6F7AEAD7DF37BAE9ECB13F679F7F7167F7E2E5F2CBA793DFE2275565BF37FAAEFFF2AEEB2E7D3EC77CFE323BF0F7646ED73F9BC37AFD5B6EBEAEB7DFC3EB1B6BFF6A51D615AE99343F7DEF667CD15103DCD13DF74D1F4C55D76B64FD14AFBD978FBCFAB9852C9D77CB9E6B9889757B4F6EB87DFB55DBFE583D9B27BBEE3AFFA79AADDFE635333D7A6EFE5AE8F9DD0B0E5D67DD26EEEAE81401FD3A9FEFE9BF38DCEAED5F4EEADFD96B917D8667938FF69D3B177ECDD6AC17FBA1BB36F8D47769F7E5F43FB793BEFA78EFFB4775F18BEFD7CFDB76E6C2FD4DB4E1C5AC30A83DEEF0498127B9798F2B7D59A259D63CFD779A9A9D85AFF3BFAF38CE6C0F6A4DEE5AB7AF1EDBED599637EF7508C2B17FE74ADC21BC2FF1BC9966BAF17F14B5FBA99D68BDE66D52FBD6E74C7B49FFD7E4ABE3F3E73DFBB6CD7BCBF5D745E169C9F6A936FD73527D8270D9388B4F5A7581F3C8E39F5AD57E95279B3B8C5BE702ADBE8F82AFB34BB1EB99BBD3B6C5B150DE9FFB9B464F4A8EA14E6DEEC475C1FABE230D9FADFA26B7FB878E226D777E6BE2F374AEB8AFE3BD3699BF182B60777B19EE7D99AAAEC6F5965836AF6784D95BDF2ADD1ED999D14E5B0ED95717D7248E6BB521FCE2A85E6B6F07C7B3A37DC14E80FCA66E8EF613FAECFAACDA9FDB5EDF9380DF1F16877D55D62F287FF3765B3E8F358C7CE47E5E4BCF0F59992B1ADF4B572EFED9EFA13EEEDE37D266F28313C15EBDFAA53E972B6EEBDB503F99A7CDB99AD9BE1E8E8E5EDDAD77E41A77D33ED89F31C7EDDF67356AB6BEDE26BA9F0BA5613DEC1FFA98BB94DBCFC38CD7F3662F73C5AAC3DE3FCF0741A2F9E2629BB4D3A8D5A33A7AAFA2E73F565ECD594BFFB08CFED71691FE45ED91D93CF1FCA8755A7D644B9D4A1F4B39D7C2C64676F70A767C69DA790FD789EE1D6E4BE6AF2BBA633373F72EDB9CC5D279537CD97737B74C2BD94507DAF55B247B30C7F7B6E3B971AE92BEA9AFC06770E47BE4DFCED8FA9D16E88B61803EC6D3FC2DEC0167E2B6D3031C6B2F3B56CEA077C935B1E9D36CFF476F86A7EE12FA6E3DFC8CD3B2E6EA8B1DE5C53FAB4DBBC4B89EBC0B90F1742832E55FFF13E3CA2D9591D1B7952E4DBAD2E924C61EE56AAF7A137B5F5FC201FB5B8F06D99795626D0F0DB7C8EE81D7439D643E4EDFD229DC2EECF2BFBD1F90339E3ABFB4361F56381F0F2E539FD4A96984B5E9620FD92CF781DBB8413FD915496E47C28FF517F77A208126FAEA39B63B8BD6CC4BF34DF07E756E2B9936DE939F88EE7E6E2E63B3B5393762BBEBDDA26EB278B26C78E875CAE7CC36E91DB96E94732EF20FF4960B5E23E67B7A1F7672FF2E1EADE3876646CCAF9DA55E6486C0BA6F0947BADD4EEE3D6C273DB95EF08DFF57AE58B856F4CEC6CF86922FFBC3B07AE37B697EDC17DD87F8B3A733AA15E9F1CF30F4BEB4CE975CEE48B969AE55691BF59CB538ED023BE8A85256D9177D225445F8F50EEAA59533FEEE21E7C817B59B392BDC4AE3F0B9E19B743D7D30EB3BED1469A0F4F459CCB0F195D571C976A3A4BC2BF1DBD53CFEE4B0BED2F7281F7D970DF3D3C1AB8D9D29E2597873E5FD609BF6927E6A1F3A3D26D2BF7235F738E9F25C7B14F7E7AE907954BD59349C3C96AC953F58AD2CB9122B077F0371D2DDF3AF730E8F2A06F52F10447C66DAC47F1C29BFE0B7389F86D6DD5A2B57145F0E5672363B221E9CFBD40D199C44649DEF895B2A55223D143709DCA85D9C43BC2A58753C35CF741FA4FC733F7E857463407B06BD0ED1E3567AE428CF3447EC5920B25174E57E5D532F9817A725B574D6AF978A90DF6609E765A4BB8E6FD38FB483BB2BDB1B159E0CF617BA385E2DBE07AAC1A1EFB65F8576C6727ABBD2548AEE764A059DB5CD955CDF29E75C6377E26CF6DA353CBDF604E09F54707B8A8D6E9FD1092F893E762AB63A43EC29C56B1DB14D3DD5024B9AF36B97ACF685F59AA51B88CB60EA97192F5C4C9FF50B4D8A4F10E59BB964C56547F2236ED26CE320F51C419F08FF9A5B8DD7AC42811B6921AA99C90F54C728FC14DC95AF64ED1BCDF40F9742763E4D574B09ABA8602DF0AA0A7DC84D4B7F1075E5F1DC1F553F38544CFD449B868B21DC497A36BF46EE354932BCF899943DB0E13955712DE8E788DAD41FD9D9EE2298B4D1AF19C35B4A8B93FB060BD7CA337AE1B1FD2C0C98FA572B1EAAA8710D918D855D15894A1B650276EA5AA3CF6AAC06E5A3DE9881379C4BE4CC7FC264886E9E2BB4C75B474D86EBBA2ADCBA9B19B186AFF27D7406F3D6F3167B347B2B9A1D8D1D0CB3B2025B76438EF178A7991CC176FB21B9EC6B0533CF55087E6D433F934DD495768E28C7B0E7B9078B1E170936450BE0C34D947D47B5DA18946B93C2CDC71E968DD634F4E16E02FEAA83B2C0B36E21578B0705D6E0DF515BE280F36855AF3015328D572DA45F259B4EE660EE9C6094FDA82DADEE9F46C468A2A4EB7A8B713BB60CED11AEE653669B570575AFEC512D6F9A42DDD2FCB7D8C948E82FE7E992EA30EE04305B361E85D54763CCEF998E84CA75FD41BF0C7A47C9423DE4AE9553DE77C75DC800713AE6092CC2EC665D1D7DE3632FE64F237BF0682F64787FD035F50D3A516999885DEE7464B1B227675D954EFF023CCB8C4E41D47DF22DA6D0133EB846A0D40BE57DB05E13F543CE574347E9FCF465D4018D7B8396239CEF19A54C6E52F3622F5481B2A2D1CAF18ED406DED50C7325B39F4238CBC0CA9DAD185F2316691AE3D0977AD1FDD7934535B990EB4983E710198D39E974483F2638047DD59D1335C698C79D57B1051E62760BA46C3435F3437E83E43A75D6DB70AB3FC1E4C28D63E3B183B163905F535C4317295E683B9CBB3591995670895A8CD1B4D67F130749AE069E727701256073BBEA0A9CA467D4470F63961BFBCC32E240B4696D09EB2901F5F6845DC37460B1AACAEC40597DB9E0737A348E2FDC03D283C38339097832B93B351659005CD398824E36FF4D7D89896CA5A98C0439E48D69246D2D67C561523D58B1403E9C119179D2754E64487C6CDC1E8886761AC857C906A9BD8BC09EDA1D4B3C8861E45D65086DF84476BE44CA2476091FD6E84A9C6A5832D44B86351189603BCF2E0EDAB5BC32DF2C717252DB6099B44AFF45B67604702F1A08C7E2D9372AD61796C867DE5AC0DB805B4606AC48E4AA647A645AF82C5F4A4B974A661CE28342A1FC5F9D0758CD01A9A72FB29E9276C2C5BCBD234FC285548446AEC990711E3CDC28C49D439515D98CEA4D222E3E62B8890AE7AC73504DD244693DA6F4B8B80E7A1E610BF129BCCADDCF51B60425071BF2CF826C1CE54FA97051B56376061892D1778F24E8319E5A87839ACC5BB0C7544F4D9306EEA097FCE66C263D1868D909599E1D598C2E826B4C51B0834B426869AB6A2F63C16AD4BFA33F61F15157F072342DB0974C0ACC9CD35320B0C8B0D48D1F2A8BCCE61A454F9F4A253DAC23E8C8830B9739468191878A670B11554D5C99CE130BC11D5C81CCD16841985B530FC49555A1F0CE365F0C9C944D585C10CA86470DC98D6871A19E83CE0E3EB65D7F8745DB8CA3BE2E3A8A3D81C5E4FCA194E4753835D685C472DE0B8E2C3FA19C00DB554AE2C8FAF2940C611791A90455EEB0112D38E631621A2B9184EA42B7682EE916A9860981D9C3F83365B9AD4F415CA4E0C5956D155F8AB94253BE095F09C78A3B9109F2B40AA8A66B0A4CF944EF5600CD34A626F900DC13552977E0757A032FDE0C7901B1C70D5BB30AAC40CDCBF22BD501C72278E722E2ACC0FC6D1634B818596DC8B4A13CC154E2719840565DA00A5D6C6840FA80A23CB4B9C1431BC404FE52C76274B24A237CC64DE499DE53B429B0A1C15938514664358996F78C08BE8081E3E88E948BB142CF19216B6BB0106ED9D0DAA367042A682AE186F2C22DF8A21F41945E9A269FC33171B816F5877A0E6411C8CB8148F031163046CA4A8DF11505D45C601AF0505AD2A8FA1D0154A10228B83487962BF6837909354E0744BC8722A5628F8912C7F010BD4FC67E05AE38209A569A8AA351A2717C0C1C62389B9055A011A0062698908DC33443B2A0E534BCBBB98D1960A932936245E3C22E110BC9148E12AD32405DF01D879A67A41B3C4915190F74B66031F64D7A985A4E91816C6072475AA8072889ECE41311BC1098DA8DE910B0B6219A9F6EC71266906594AD68A14CBB29A86D28AF1422660A5651C7C3BBBCC853020AAE0329A382C009760A4EE95047C097BDC6134E29D6A68BA3473A04B945688982D0D4C03869F10B52AEAE10309822EA56D658C662A526B176D107033C70252CD42CAE0E119434A5398493E02B8CF0452C02BC3405E02746A84B75E18359D58CB41D188CF5043E0D8652A41A05C0996028214415A1873496C70CFD0E4328011D162CE3DC1124D47AC8164193E452A017C6894581E48634AF00412E00563A996D7929A1C3E5A06E822C4C3C0B0A022D059CAC107806118C8ADC86835C282B0E856D3C0A8F109D3F30640D51424EC44594407841F480D091E8243B7804DB94631C54F584B282550D114A6E29104EB37F424A84C02ED29C4258929F19665463932A496297A14BF1AFF8E9AA46566C5D7A2FD34CC0A3B68393109B803C61604881D4914525A9242A393A668649A17B0478512030268C815D85700D10AC8894875134498A0A5BCE88CC1A9484919C5D584416A3C07FC85534D82CCC2842D8454C451864B9062B79190B153ABC1E222733204B5009CEA2289909160C488AD90491378D6C140DA439804B210D74213C5295DE03587EA4C4CE8298D2062AC6220693440AD512858FCCB70FD180538648C99E0490E05166C07C610F2186748422BA9131DD481460C29A072E4803AD80C8A6D60FBD00F298FD70148113E18B094C18EC205D9271664229EA20DA09F149848E0DD5FB2B1183F09293F8201A8C5EA329964BF9224150211A65CCFA0A2794993BC45D001231E2744030546148DD374504F5469561B52D04B0864D4965AD62D770015148598301C2E4C04D44E841201D08A0C3B5800F99489287942E54A02D656C4D51110A8095A3B4E7373934036493842C8D2810769139710874A1D321913419AD108FD6A21C61C9C8A39014D506C6A1470C41AD70612BD218A61C128736B4E02D209E031490509B6254D6040B1D030438ACC47B44925058904A4A380F5AC99B207E820C43FD8C57D250080E9C4AAF4C23D48C60D984FA6A4502A402E286A8022B3BB090B981B16A1004A21901CD08D8373875064970B0801EF843856053261B8B0187DC320A514FD022105CB3D30731482E1C302A8853E0E00D8210AF22504109074D09C92BCB38B601F04402632927A9EC49A78D2C8800CB16088A18481DAA93004F401C8404852199241CFF280C401A72043DA086A1389A49BC0A1EA693223082359A0808893941138118AA6081C0466E20BD12AE71A4A7091B9AA890D30D45}


\providecommand{\bysame}{\leavevmode\hbox to3em{\hrulefill}\thinspace}
\providecommand{\MR}{\relax\ifhmode\unskip\space\fi MR }
\providecommand{\MRhref}[2]{%
  \href{http://www.ams.org/mathscinet-getitem?mr=#1}{#2}
} \providecommand{\href}[2]{#2}

\end{document}